\documentclass[12pt]{amsart}
\usepackage{fullpage}
\usepackage{graphicx}
\usepackage{caption, subcaption}
\usepackage{amssymb}
\usepackage{epstopdf}
\usepackage{float}
\usepackage{enumerate}

\def\qed{\hfill
\ifhmode\unskip\nobreak\fi\quad\ifmmode\Box\else$\Box$\fi\\ }
\newtheorem{thm}{Theorem}
\newtheorem{cor}[thm]{Corollary}
\newtheorem{lem}[thm]{Lemma}

\begin{document}

\author{Ervin Gy\H{o}ri}
\address{Alfr\'{e}d R\'{e}nyi Institute of Mathematics \\ Budapest, Hungary \\ and\\ Department of Mathematics, Central European University \\ Budapest, Hungary}
\thanks{Research of this author is supported in part by OTKA Grants 78439 and 101536}
\email{ervin@renyi.hu}
\author{Alexandr Kostochka}
\address{Department of Mathematics \\ University of Illinois \\ Urbana, IL 61801, USA\\and\\ Sobolev Institute of Mathematics\\ Novosibirsk, Russia}
\thanks{Research of this author is supported in part by NSF grant  DMS-1266016 and  by grants 12-01-00631 and 12-01-00448 of the Russian Foundation for Basic Research. }
\email{kostochk@math.uiuc.edu}
\author{Andrew McConvey}
\address{Department of Mathematics \\ University of Illinois \\ Urbana, IL 61801\\ USA}
\email[Corresponding author]{mcconve2@illinois.edu}
\author{Derrek Yager}
\address{Department of Mathematics \\ University of Illinois \\ Urbana, IL 61801\\ USA}
\thanks{The author acknowledges support from National Science Foundation grant DMS 08-38434 ``EMSW21-MCTP: Research Experience for Graduate Students.''}
\email{yager2@illinois.edu}
\title{A list version of graph packing}
\begin{abstract} We consider the following generalization of graph packing.
Let $G_{1} = (V_{1}, E_{1})$ and $G_{2} = (V_{2}, E_{2})$ be graphs of order $n$ and $G_{3} = (V_{1}\cup V_{2},E_{3})$ a bipartite graph.  A bijection $f$ from $V_{1}$ onto $V_{2}$ is a \emph{list packing} of the triple $(G_{1}, G_{2}, G_{3})$ if $uv \in E_{2}$ implies $f(u)f(v) \notin E_{2}$ and $vf(v) \notin E_{3}$ for all $v \in V_{1}$.  
We extend the classical results of Sauer and Spencer and Bollob\' as and Eldridge on packing of graphs with small sizes or maximum degrees to the
setting of list packing. In particular, we extend the well-known Bollob{\'a}s--Eldridge Theorem, proving that if $\Delta (G_{1}) \leq n-2, \Delta(G_{2}) \leq n-2, \Delta(G_{3}) \leq n-1$, and $|E_1| + |E_2| + |E_3| \leq 2n-3$, then either $(G_{1}, G_{2}, G_{3})$ packs or is one of 7 possible exceptions. Hopefully, the concept of list packing will help
to solve some problems on ordinary graph packing, as the concept of list coloring did for ordinary coloring.
\end{abstract}
\maketitle

{\small{Mathematics Subject Classification: 05C70, 05C35.}}{\small \par}

{\small{Keywords: Graph packing, maximum degree, edge sum, list coloring.}}{\small \par}

%
%
\section{Introduction}\label{sec:intro}
%
%
The notion of graph packing is a well-known concept in graph theory and combinatorics.  
Two graphs on $n$ vertices are said to \emph{pack} if there is an edge-disjoint placement of the graphs onto the same set of vertices.  
In 1978, two seminal papers,~\cite{S-S} and~\cite{B-E}, on extremal problems on graph packing appeared in the same journal.
In particular, Sauer and Spencer \cite{S-S} proved sufficient conditions for packing two graphs with bounded product of maximum degrees.  
\begin{thm}[\cite{S-S}]\label{S-S product}
Let $G_1$ and $G_2$ be two graphs of order $n$. If $2\Delta(G_1)\Delta(G_2)<n,$ then $G_1$ and $G_2$ pack.
\end{thm}

This result is sharp and later Kaul and Kostochka \cite{K-K}  characterized all graphs in which Theorem~\ref{S-S product} is sharp.

\begin{thm}[\cite{K-K}]\label{K-K}
Let $2\Delta(G_1)\Delta(G_2)\leq n.$ $G_1$ and $G_2$ do not pack if and only if one of $G_1$ and $G_2$ is a perfect matching and the other is either $K_{\frac{n}{2},\frac{n}{2}}$ with $\frac{n}{2}$ odd or contains $K_{\frac{n}{2}+1}.$
\end{thm}

In the same paper, Sauer and Spencer  gave sufficient conditions for packing two graphs with given total number of edges. 

\begin{thm}[\cite{S-S}]\label{S-S}
Let $G_1$ and $G_2$ be two graphs of order $n.$ If $|E(G_1)| + |E(G_2)| \leq \frac{3}{2}n -2$, then $G_1$ and $G_2$ pack.
\end{thm}

This result is best possible, since $G_{1} = K_{1, n-1}$ and $G_{2} = \frac{n}{2} K_{2}$ do not pack.  
Independently,  Bollob{\'a}s and Eldridge \cite{B-E} proved the stronger result that the bound of Theorem~\ref{S-S} can be significantly strengthened when $\Delta(G_1) < n-1$
and $\Delta(G_2) < n-1$. 

\begin{figure}[!ht]
\begin{subfigure}[t]{.22\textwidth}
  \centering
  \includegraphics[width=.9\linewidth]{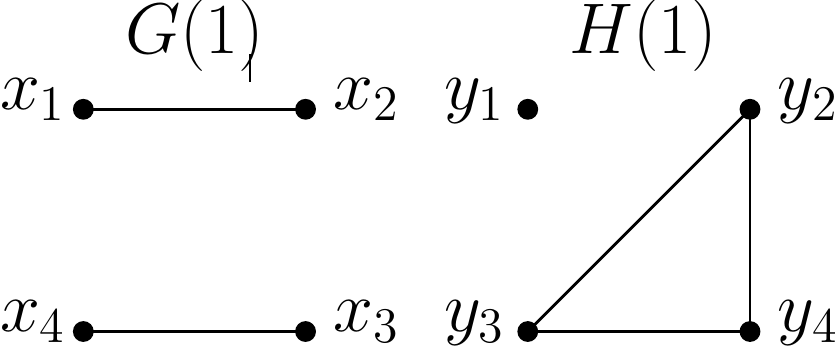}
  \label{fig:BE1}
\end{subfigure}%
\begin{subfigure}[t]{.22\textwidth}
  \centering
  \includegraphics[width=.9\linewidth]{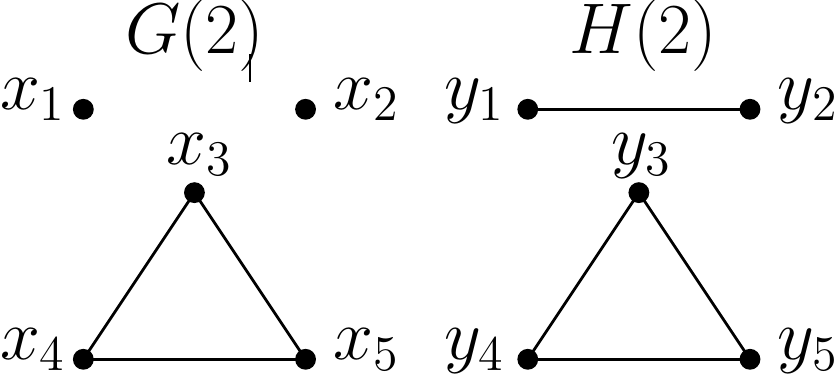}
  \label{fig:BE2}
\end{subfigure}
\begin{subfigure}{.22\textwidth}
  \centering
  \includegraphics[width=.9\linewidth]{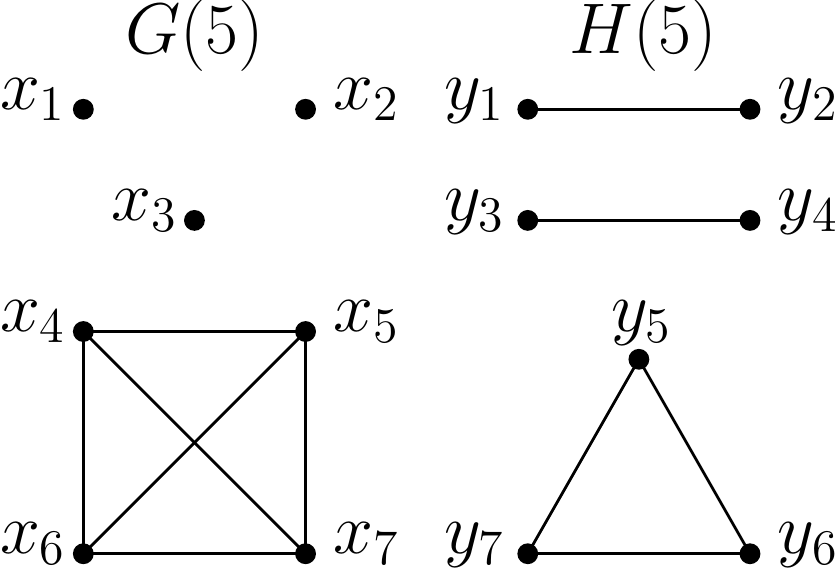}
  \label{fig:BE5}
\end{subfigure}
\begin{subfigure}{.22\textwidth}
  \centering
  \includegraphics[width=.9\linewidth]{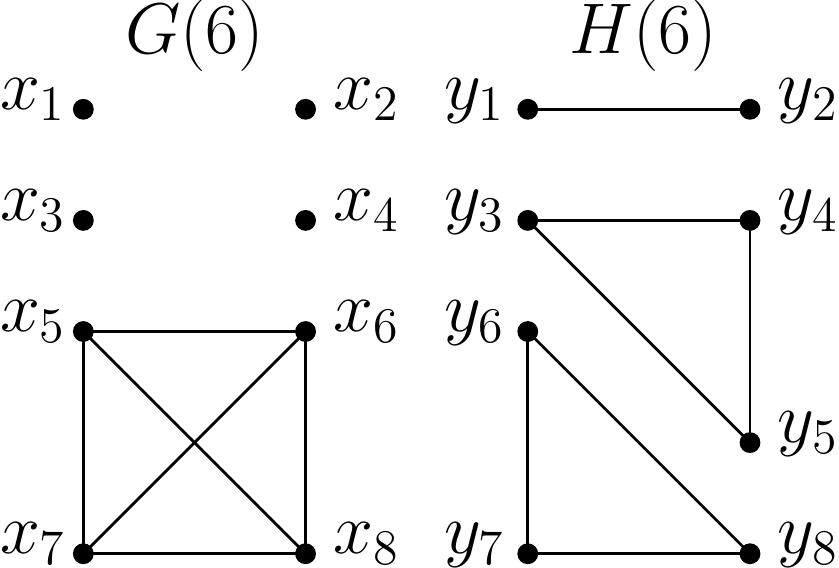}
  \label{fig:BE6}
\end{subfigure}

\begin{subfigure}[b]{.22\textwidth}
  \centering
  \includegraphics[width=.9\linewidth]{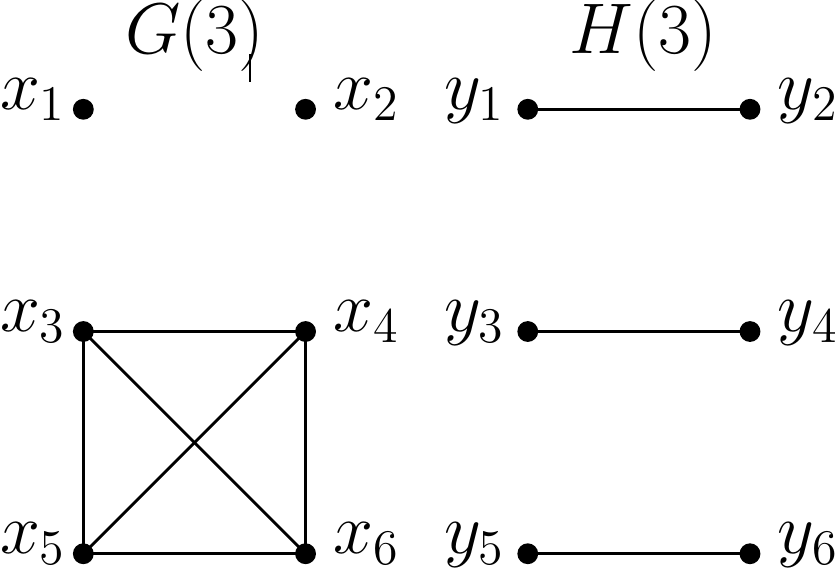}
  \label{fig:BE3}
\end{subfigure}
\begin{subfigure}[b]{.22\textwidth}
  \centering
  \includegraphics[width=.9\linewidth]{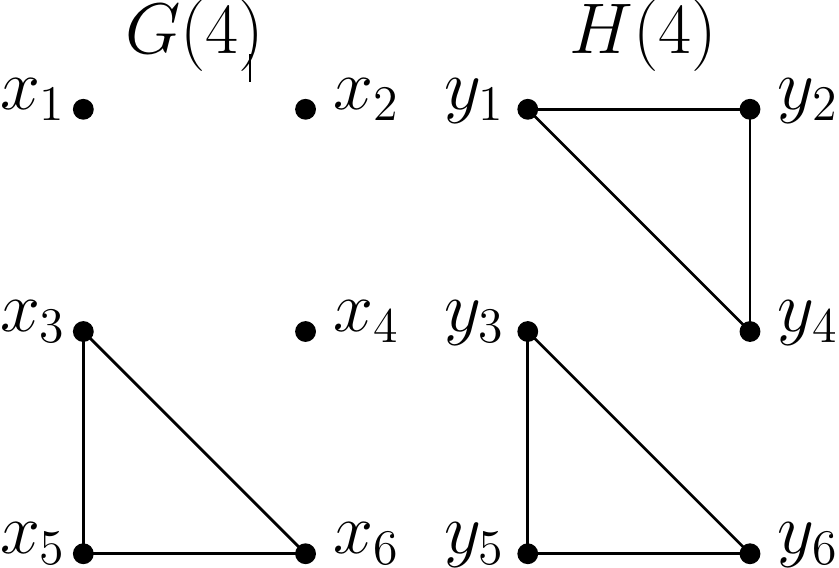}
  \label{fig:BE4}
\end{subfigure}
\begin{subfigure}[b]{.44\textwidth}
  \centering
  \includegraphics[width=.9\linewidth]{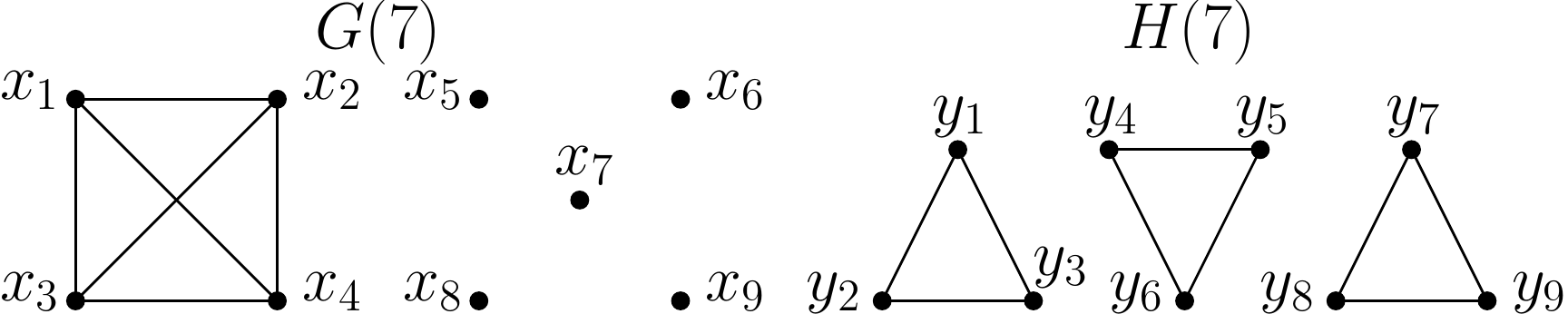}
  \label{fig:BE7}
\end{subfigure}
\caption{Bad pairs in Theorems~\ref{B-E} and~\ref{List B-E}.}
\label{fig:B-E}
\end{figure}

\begin{thm}[\cite{B-E}]\label{B-E}
If $\Delta(G_1),\Delta(G_2) \leq n-2, e(G_1)+e(G_2)\leq 2n-3,$ and $\{G_1,G_2\}$ is not one of the following pairs: $\{2K_2,K_1\bigcup K_3\}, \{\overline{K_2}\bigcup K_3, K_2\bigcup K_3\}, \{3K_2, \overline{K_2}\bigcup K_4\}, \{\overline{K_3}\bigcup K_3,2K_3\},$ $\{2K_2\bigcup K_3,\overline{K_3}\bigcup K_4\},\{\overline{K_4}\bigcup K_4, K_2\bigcup2K_3\},\{\overline{K_5}\bigcup K_4, 3K_3\}$ (Figure~\ref{fig:B-E}). Then, $G_1$ and $G_2$ pack.
\end{thm}
This result is also sharp, since the graphs $G_{1} = C_{n}$ and $G_{2} = K_{1,n-2} \cup K_{1}$ satisfy the maximum degree condition, have $2n-2$ edges, and do not pack.
There are other extremal examples.

Variants of the packing problem have been studied and, in particular, restrictions of permissible packings arise both within proofs and are posed as independent questions.  The notion of a bipartite packing was introduced by Catlin \cite{C} and was later studied by Hajnal and Szegedy \cite{H-S}.  This variation of traditional packing involves two bipartite graphs $G_{1} = (X_{1} \cup Y_{1}, E_{1})$ and $G_{2} = (X_{2} \cup Y_{2}, E_{2})$ where permissible packings send $X_{1}$ onto $X_{2}$ and $Y_{1}$ onto $Y_{2}$.  The problem of fixed-point free embeddings, studied by Schuster in 1978, considers a different restriction to the original packing problem \cite{Schuster}.  In this case, two edge disjoint copies of a graph $G$ are placed into $K_{n}$ with the additional property that two copies of the same vertex must be mapped to different vertices in $K_{n}$.  In \cite{Z}, Schuster's result is used to prove a necessary condition for packing two graphs with given maximum and average degrees.

In this paper, we introduce the language of list packing in order to model such problems.  A \emph{list packing} of the graph triple $(G_1, G_2, G_3)$ with $G_1=(V_1,E_1), G_2=(V_2,E_2),$ and $G_3=(V_1\bigcup V_2,E_3)$ is a bijection $f: V_{1} \rightarrow V_{2}$ such that $uv \in E_{1}$ implies $f(u)f(v) \notin E_{2}$ and for each $u \in V_{1}$,  $uf(u)\notin E_3$.  Note that both $G_1$ and $G_2$ are graphs on $n$ vertices so that $G_3$ has $2n$ vertices, and one can think of the edge set $E_{3}$ as a list of restrictions that must be avoided when packing $G_1$ and $G_2$.

This notion is closely related to Vizing's concept of list coloring \cite{V}.  Suppose we wish to color a graph $G$ with the colors $\{1, \ldots, k \}$.  A \emph{list assignment} $L$ is a function on the vertex set $V(G)$ that returns a set of colors $L(v) \subset \{1, \ldots, k \}$ \emph{not} permissible for $v$.  A \emph{list coloring}, more specifically an \emph{$L$-coloring}, is a proper coloring $f$ of $G$ such that $f(v) \notin L(v)$ for all $v\in V(G)$.  In fact, the problem of list coloring $G$ can be stated within the framework of list packing.  A proper $L$-coloring of a graph $G$ is equivalent to a list packing where $G_1=G$ along with an appropriate number of isolated vertices, $G_2$ is a disjoint union of $K_{n}$'s each representing a color, and $E_3$ consists of all edges going between a vertex $v \in V_1$ and the copies of $K_{n}$ corresponding to colors in $L(v)$.

Similarly, the restrictions to packings discussed above can be modeled using this framework.  A bipartite packing is a packing of the triple $(G_{1}, G_{2}, G_{3})$ where $E_{3}$  consists of all edges  between $X_{i}$ and $Y_{3-i}$ for $i = 1, 2$.  A fixed-point free embedding is a packing of the triple $(G, G, G_{3})$ where $E_{3} = \{ (v,v): v \in V(G)\}$.  
Although the list packing is more general, some important theorems on the ordinary packing can be transferred to the list setting.
The results of this paper prove natural generalizations of Theorems~\ref{S-S product}--\ref{B-E} in the language of list packing.  In particular, we extend Theorem~\ref{S-S product} and Theorem~\ref{K-K} as follows.

\begin{thm}\label{List S-S product}
 Let $G = (G_{1}, G_{2}, G_{3})$ be a graph triple with $|V_{1}| = |V_{2}| = n$.  If $\Delta (G_{1}) \Delta (G_{2}) + \Delta (G_{3}) \leq n/2,$ then $G$ does not pack if and only if $\Delta (G_3) = 0$ and one of $G_1$ or $G_2$ is a perfect matching and the other is $K_{\frac{n}{2},\frac{n}{2}}$ with $\frac{n}{2}$ odd or contains $K_{\frac{n}{2}+1}.$ Consequently, if $\Delta(G_{1}) \Delta (G_{2}) + \Delta (G_{3}) < n/2,$ then $G$ packs.\end{thm}

The main result of this paper is the following list version of Theorem~\ref{B-E}.

\begin{thm}\label{List B-E} Let $n\geq 1$ and $G_1$ and $G_2$ be $n$-vertex graphs.  If $\Delta (G_1),\Delta (G_2) \leq n-2$, $\Delta (G_3) \leq n-1$, $|E(G_1)| + |E(G_2)|+|E(G_3)| \leq 2n-3$ and the pair $(G_1,G_2)$ is none of the 7 pairs in Figure~\ref{fig:B-E}, then $G_1$ and $G_2$ pack.
\end{thm}

\begin{figure}[h!]
\begin{subfigure}{.33\textwidth}
  \centering
  \includegraphics[width=.5\linewidth]{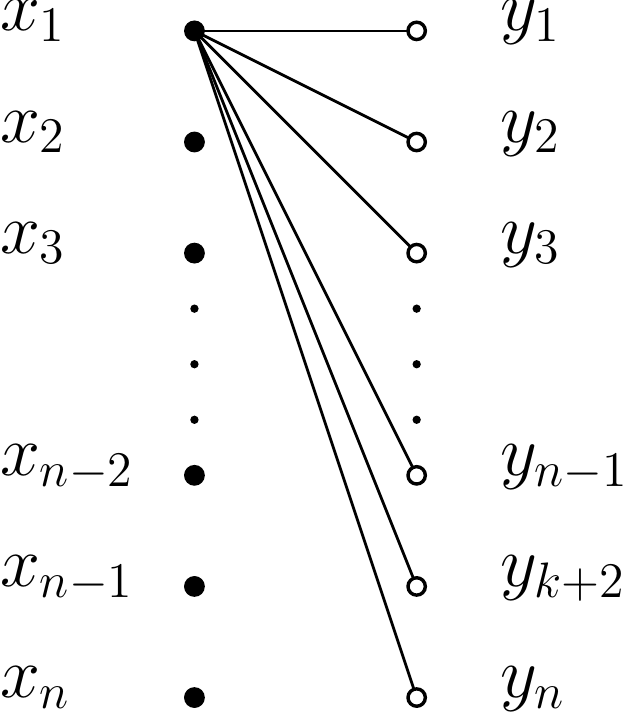}
  \label{fig:yellow star}
  \caption{}
\end{subfigure}%
\begin{subfigure}{.33\textwidth}
  \centering
  \includegraphics[width=.5\linewidth]{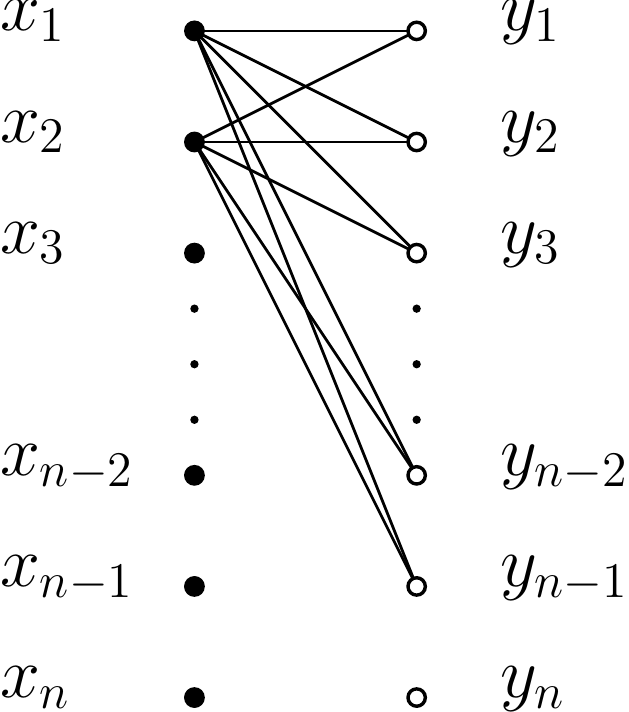}
  \label{fig:double yellow star}
  \caption{}
\end{subfigure}
\begin{subfigure}{.33\textwidth}
  \centering
  \includegraphics[width=.5\linewidth]{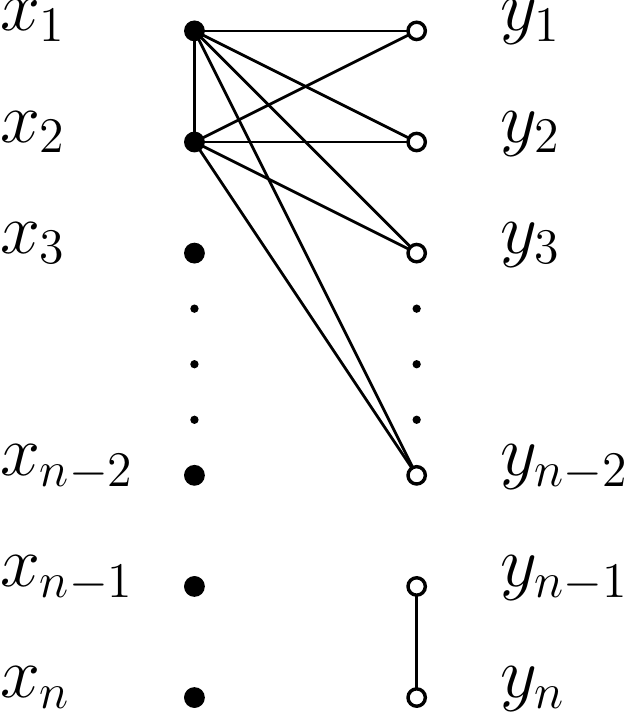}
  \label{fig:two edges}
  \caption{}
\end{subfigure}
\begin{subfigure}{.33\textwidth}
  \centering
  \includegraphics[width=.5\linewidth]{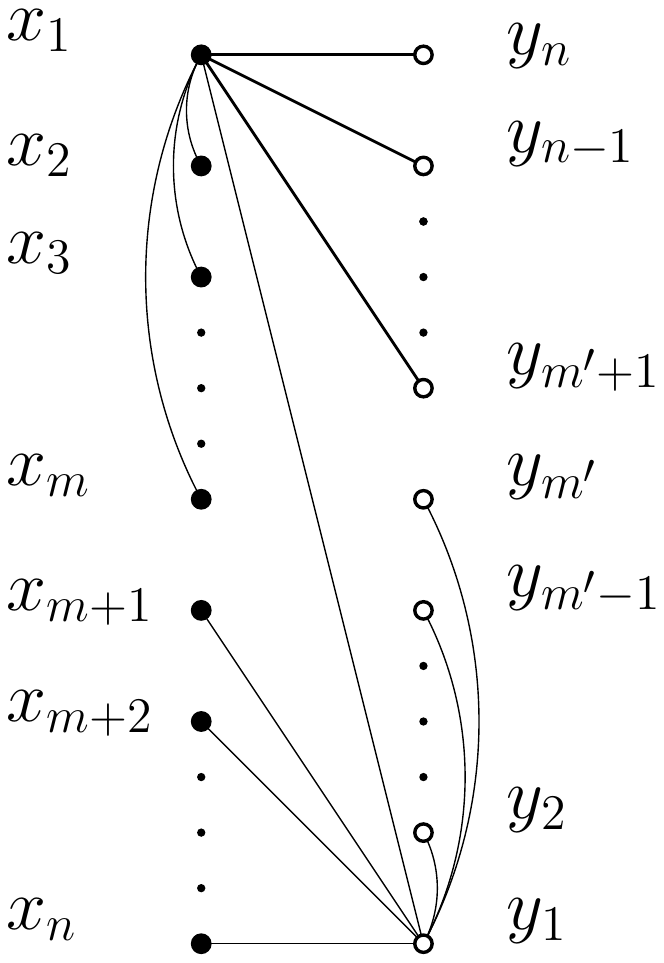}
  \label{fig:partition}
  \caption{}
\end{subfigure}
\begin{subfigure}{.33\textwidth}
  \centering
  \includegraphics[width=.5\linewidth]{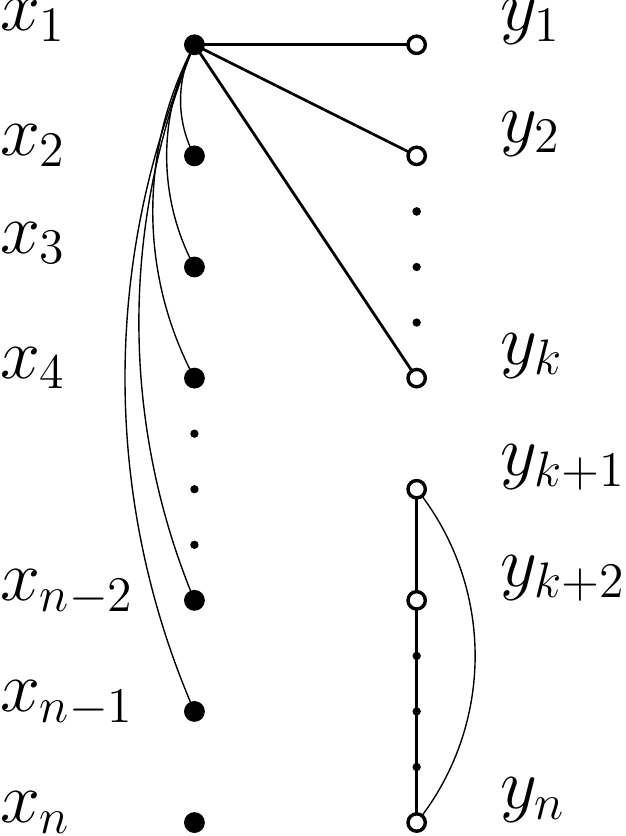}
  \label{fig:cycle}
  \caption{}
\end{subfigure}
\caption{Sharpness examples for Theorem~\ref{List B-E}}
\label{fig:sharpness examples}
\end{figure}

Theorem~\ref{List B-E} is sharp and the list version introduces several new sharpness examples.  First, the condition $\Delta_{3} \leq n-1$ cannot be removed, since a vertex in $V_{1}$ adjacent to all vertices in $V_{2}$ cannot be placed at all (Figure~\ref{fig:sharpness examples}A).  The restriction on the edge sum is also sharp, as there are  several examples of graphs with $E_{3} >0$ and edge sum equal to $2n-2$ that do not pack.  We provide 4 such examples.

For the first example, consider $G_{1}$ and $G_{2}$ to be independent sets and $x_1, x_2 \in V_{1}$ each adjacent to all but a single vertex of $V_{2}$ (Figure~\ref{fig:sharpness examples}B).  Alternatively, consider $E_{1}$ consisting of a single edge $x_{1}x_{2}$, $E_{2}$ consisting of a single edge $y_{n-1}y_{n}$, and $E_{3}$ consisting of all edges between $\{ x_{1}, x_{2} \}$ and $V_{2} - y_{n-1} - y_{n}$ (Figure~\ref{fig:sharpness examples}C).  For the third example, consider  $G_{1} = K_{1, m-1} \cup \overline{K}_{n-m}$, $G_{2} = K_{1, m'-1} \cup \overline{K}_{n-m'}$ (for any choice of $m,m'$), and $E_{3}$ consisting of all yellow edges between the center of the star in $G_{1}$ and isolated vertices in $V_{2}$ as well as between the center of the star in $G_{2}$ and isolated vertices in $V_{1}$ (Figure~\ref{fig:sharpness examples}D).  Finally, consider $G_{1} = K_{1, n-1} \cup K_{1}$, $G_{2} = C_{k} \cup \overline{K}_{n-k}$ (for any choice of $k$), and let $E_{3}$ consist of all possible edges between the center of the star in $G_{1}$ and isolated vertices in $G_{2}$ (Figure~\ref{fig:sharpness examples}E).

Though this paper  focuses on extending classical packing results to the list setting, one of our goals is to provide tools to handle
 problems of standard graph packings.  In particular, we heavily use Theorems~\ref{List S-S product} and~\ref{List B-E} in~\cite{ListZak}
   to get an approximate solution to a conjecture of \.{Z}ak~\cite{Z} on packing $n$-vertex graphs with given sizes and maximum degrees.

The paper is organized as follows.  In the next paragraph, we introduce some notation.
In Section~\ref{sec:List S-S product}, we prove Theorem~\ref{List S-S product}.  Section~\ref{sec:Prelim} contains some preliminary results, including an extension of Theorem~\ref{S-S} that will be used as a base case in our proof of Theorem~\ref{List B-E}.  Then, Section~\ref{sec:Proof of list B-E} contains our proof of the main result  by induction on the size of the vertex set.

%
\subsection{Notation}\label{sec:notation}
%
A graph triple $G = (G_{1}, G_{2}, G_{2})$ of size $n$ consists of a pair of $n$-vertex graphs $G_{1} = (V_1, E_1)$ and $G_{2} = (V_{1}, E_{2})$ together with a bipartite graph $G_{3} = (V_{1} \cup V_{2}, E_{3})$.  Let $V:= V_{1} \cup V_{2}$.  An edge in $E_{1} \cup E_{2}$ is a \emph{white} edge, while an edge in $E_{3}$ is a \emph{yellow} edge.  For $v \in V_{i}$ ($i = 1,2$), the \emph{white neighborhood} of $v$, denoted $N_{i}(v) \subseteq V_{i}$, is the set of neighbors of $v$ in $G_{i}$, $d_{i}(v) = |N_{i} (v)|$, and $\Delta_i = \max_{v\in V_i}d_i(v)$.  For convenience, when $w \in V_{3-i}$, we say that $N_{i}(w) = \emptyset$ (and hence $d_{i}(w) = 0$).  The \emph{yellow neighborhood} of $v$, denoted $N_{3}(v) \subseteq V_{3-i}$ is the set of neighbors of $v$ in $G_{3}$ and $d_{3} (v) = |N_{3}(v)|$.  For $v \in V_{i}$, the \emph{neighborhood} in $v$, denoted $N(v)$ is the disjoint union $N_{i}(v) + N_{3}(v)$ and the \emph{degree} of $v$ is $d_{i}(v) + d_{3}(v)$ and is denoted $d(v)$.

For $i = 1,2,3$, let $e_{i} = |E_{i}|$ and define $\Delta_{i}$ to be the $\max_{v \in V} d_{i}(v)$.  Finally, the triple $G$ \emph{packs} if there is a bijection $f:V_{1} \rightarrow V_{2}$ such that $v f(v) \notin E_{3}$ for any $v \in V_{1}$ and $uv \in E_{1}$ implies $f(u)f(v) \notin E_{2}$.

%
%
\section{Proof of Theorem~\ref{List S-S product}}\label{sec:List S-S product}
%
%
\noindent$(\Leftarrow)$ Suppose $G_1$ is a perfect matching.  If $G_2$ contains $K_{\frac{n}{2}+1},$ then for any mapping $f:V_1\rightarrow V_2,$ some edge of $G_1$ will be mapped to an edge in the clique.  Otherwise, $G_2$ is $K_{\frac{n}{2},\frac{n}{2}}$ with $\frac{n}{2}$ odd, then under any mapping, we are again forced to have some matching edge in $G_1$ mapped so that it has one endpoint in each partite set. 

\vspace{.1in}

\noindent$(\Rightarrow)$ Assume that our graph triple $G$ is the minimal counterexample that does not pack where we interpret minimal as the minimal number of total edges.  If $\Delta_3=0,$ then the result follows from Theorem~\ref{K-K}.  Hence, we can assume $E_3\neq\emptyset$. By minimality, we may assume that there is a partial packing $f$ which has a conflict at only a single edge in $vw \in E_3$, where $f(v)=w$. For an arbitrary $a \in V_1-v$ with $f(a) = b$, define the mapping $f_a$ by $f_a(v)=b,f_a(a)=w$ and $f_a=f$ otherwise.  In particular, $f_a$ will be a packing of the graph triple $G$ if $a$ satisfies:
\begin{enumerate}[(i)]
\item $f_{a}(N_1 (a)) \cap N_2 (w) = \emptyset$,
\item $f_{a}(N_1 (v)) \cap N_2 (b) = \emptyset$,
\item $b \notin N_3(v),$ and
\item $w \notin N_3(a)$
\end{enumerate}
Note that there are at most $\Delta_1 \Delta_2$ vertices in $V_1-v$ that may violate (i) [similarly for (ii)] and at most $\Delta_3-1$ vertices in $V_1-v$ that may violate (iii) [similarly for (iv)].  Since $G$ does not pack, $(n-1)-[(\Delta_3-1)+(\Delta_3-1)+2\Delta_1\Delta_2]\leq0.$  But this inequality yields $n+1\leq2[\Delta_3+\Delta_1\Delta_2],$ a contradiction.\qed

%
%
\section{Preliminary facts}\label{sec:Prelim}
%
%
The following lemma is an extension of Theorem~\ref{S-S}.

\begin{lem}\label{List S-S}
Let $G_{1} = (V_{1}, E_{1})$ and $G_{2} = (V_{2}, E_{2})$ be graphs of order $n$ and let $G_{3} = (V_{1} \cup V_{2}, E_{3})$ be a bipartite graph with partite sets $V_{1}$ and $V_{2}$.  If $d_{3}(v) \leq n-1$ for each $v \in V_{1} \cup V_{2}$ and $e_{1} + e_{2} + e_{3} \leq \left\lfloor\frac{3}{2}n\right\rfloor -2$,  then the triple $G=(G_{1},G_{2},G_3)$ packs.
\end{lem}

\noindent\textbf{Proof:}
If $e_{3} = 0$, then the result holds from Theorem~\ref{S-S}.  Further, if $e_{i} = 0$ for $i \in \{1,2\}$, then the problem reduces to finding a matching in $G_{3}$ which can be done by Hall's Theorem.  So we assume that $e_{1}, e_{2}, e_{3} > 0$.

\begin{figure}[h]
\begin{center}
\includegraphics[scale=.75]{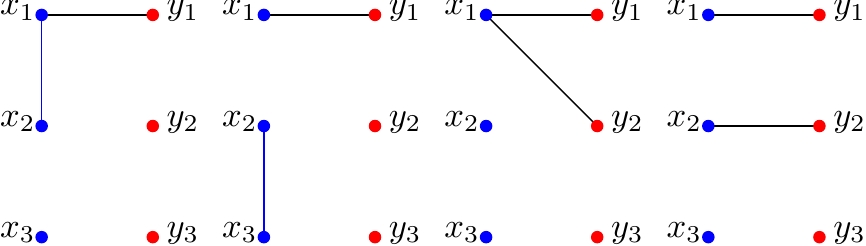}
\end{center}
\caption{Graphs with $n=3$ vertices and $2$ edges.}
\label{fig:SSbasecase}%
\end{figure}

It is sufficient to prove the case when $e_{1} + e_{2} + e_{3} = \left\lfloor\frac{3}{2}n\right\rfloor -2$.    The proof will proceed by induction on $n$.  If $n =2$, then $e_{1} + e_{2} + e_{3} = 1$ and it is clear that there is a packing.  Similarly, if $n=3$, then $e_{1} + e_{2} + e_{3} = 2$ and, up to isomorphism, there are $4$ cases.  It can be easily checked that there is a packing in each of these cases (Figure~\ref{fig:SSbasecase}).

Now assume that $n \geq 4$ and the theorem holds for all smaller values.  If there is some $v \in V_{i}$ with $d_{3} (v) = n-1$, then there are fewer than $n/2 -1$ edges not adjacent to $v$.  Let $u \in V_{3-i}$ be a vertex who has no neighbors in $(V_{1} \cup V_{2}) - v$.  If $uv \notin E_{3}$, then we pack $G_{i} - v$ and $G_{3-i} - u$ by induction and extend it by assigning $v$ to $u$.  If $uv \in E_{3}$, then there is some vertex $w \in V_{i} - v$ with degree at least $1$, otherwise we could easily send $v$ to its non-neighbor and $V_1-v$ can be sent arbitrarily. So, with this $w\in V_i-v$ where $d(w)\geq1$, we assign $w$ to $u$ and then pack $G_{i} - w$ and $G_{3-i} - u$ by induction.  We can now assume that for all $v \in V_{1} \cup V_{2}$, $d_{3}(v) \leq n-2$ and proceed in three cases:

\textbf{Case 1:} \emph{There exists a vertex $x \in V_{i}$ such that $d_{i}(x) = d_{3}(x) = 0$.} If there exists a $y \in V_{3-i}$ with $d_{3-i}(y) + d_{3}(y) \geq 2$, then $G_{i} - x$ and $G_{3-i}-y$ pack by induction and this packing can be extended to the original graphs by assigning $x$ to $y$.  So we may assume that $d_{3-i}(y) + d_{3}(y)\leq 1$ for all $y \in V_{3-i}$ and we can assume there is at least one $y\in V_{3-i}$ with $d_{3-i}(y) + d_{3}(y)= 1$ or else the graph triple packs trivially.  Moreover, we can assume there is at least one $y\in V_{3-i}$ with $d_{3-i}(y)=1$, say $yz\in E_{3-i}$ or else all edges incident to $V_{3-i}$ vertices are yellow and this too packs trivially since $G_3$ would be a yellow matching with $E_2=\emptyset$ . Since at most $n$ edges are accounted for with endpoints in $V_{3-i},$ then there is a $w\in V_{i}$ with $d_{i}(w)+d_{3}(w)\geq 2$ and $wz \notin E_{3}$.  Then, by induction, there is a packing $G_{i} - \{w,x\}$ and $G_{3-i} - \{ y,z\}$ which can be extended to a packing of the original graphs by mapping $w$ to $z$ and $x$ to $y$.

\textbf{Case 2:} \emph{There is some $x \in V_{i}$ with $d_{i}(x) = 0$, but $d_{3}(x) > 0$.} If $d_{3}(x) \geq 2$, then we find an allowed partner $z \in V_{3-i}$, pack $G_{i} -x$ and $G_{3-i} -z$ by induction, and extend the packing by assigning $x$ to $z$.  So, we may assume $d_{3}(x) = 1$.  Let $xy \in E_{3}$ be this edge.  Since there are no isolated vertices in $V_{3-i}$ (as otherwise we would be done by Case 1), let $v \in V_{3-i} - y$ such that $d_{3-i}(v) + d_{3}(v) \geq 1$.  Then, there is a packing of $G_{i} - x$ and $G_{3-i} - v$ by induction which can be extended to the original graphs by mapping $x$ to $v$.

\textbf{Case 3:} \emph{$\delta_{1} > 0$ and $\delta_{2} > 0$.} Without loss of generality, assume that $e_{1} \leq e_{2}$, so $e_{1} < 3n/4$.  There are more than $n/4$ non-trivial tree components in $G_1$ and, since $\delta_{1} \geq 1$, more than $n/2$ vertices of degree $1$.  Further, $e_{3} \leq \left\lfloor\frac{3}{2}n\right\rfloor -2 - e_{1} - e_{2} < n/2$, so there exists a vertex $x\in V_1$ with $d(x)=1$ so $x$ only has a white neighbor. Let $y \in G_{2}$ such that $d_{3}(y) \geq 1$.  Consider the graph obtained by removing $\{x,y\}$ from $V_{1} \cup V_{2}$ and adding to $E_{3}$ all edges from $N(x)$ to $N(y)$.  This results in a net change of at least two fewer edges so that, by induction, there is a packing of $G_{1} - x$ and $G_{2} - y$ which extends to a packing of the original graphs by mapping $x$ to $y$.
\qed

Lemma~\ref{List S-S} along with the following corollary will serve as a base case for our proof of Theorem~\ref{List B-E}.

\begin{cor}\label{n edges}
Suppose $|G_{1}| = |G_{2}| = n \geq 2$ and $G$ is a triple $(G_{1}, G_{2}, G_{3})$.  If $e_{1} + e_{2} + e_{3} \leq n$, then either:

\begin{enumerate}
\item $G$ has a packing, or
\item For some $i \in \{ 1, 2\}$, some $v \in V_{i}$ is adjacent to all vertices in $V_{3-i}$, or
\item $n = 2$ and $G_{1} \cong G_{2} \cong K_{2}$.
\end{enumerate}
\end{cor}

\begin{figure}[ptbh]
\begin{center}
\includegraphics[scale=.75]{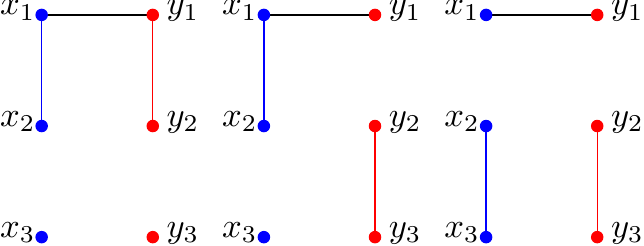}
\end{center}
\caption{Graphs with $n=3$ vertices and $3$ edges.}
\label{BEbasecase}%
\end{figure}

\noindent{\bf  Proof:}  If $n = 2$ and $e_{1} + e_{2} + e_{3} = 2$, then the result is clear.  If $n=3$ and $G_{i}$ has no white edges for some $i \in \{1,2\}$, then the problem is equivalent to finding a matching in the complement of $G_{3}$ and the result follows from Hall's Theorem.  Similarly, if there are no yellow edges, then the result follows from Figure~\ref{fig:B-E}, so it must be the case that $e_{1} = e_{2} = e_{3} = 1$.  Up to isomorphism, there are only 3 cases and it is clear that in each case there is a packing (Figure~\ref{BEbasecase}).   For $n \geq 4$, the result follows from Lemma~\ref{List S-S}.  \qed

%
%
\section{Proof of Theorem~\ref{List B-E}}\label{sec:Proof of list B-E}
%
%
Let $G=(G_1,G_2,G_3)$ of size $n$ be a counterexample to Theorem~\ref{List B-E} with the smallest size.  Notice that $n \geq 4$, since otherwise Corollary~\ref{n edges} applies. Also, we assume $E_3\neq\emptyset$ or else Theorem~\ref{B-E} applies.

\begin{lem}\label{yellow star}
$\Delta_3\leq n-2$.
\end{lem}

\noindent\textbf{Proof:}  Suppose that there exist $v \in V_1$ and $w\in V_2$ such that $N_{3}(v)=V_2-w$. Let $G'$ be obtained from $G$ by deleting all $n-1$ edges connecting $v$ with $V_2$ and all edges (maybe zero) connecting $w$ with $V_1$. Let $A=N_{G'}(v)$, $B=N_{G'}(w)$, $a=|A|$ and $b=|B|$. If $a=0$ or $b=0$, then if we pack $G_1-v$ with $G_2-w$, placing $v$ on $w$ does not create conflicts. And $G-v-w$ has at most $(2n-3)-(n-1)=n-2$ edges. Such pairs always pack. So assume $a\geq 1$ and $b\geq 1$.

Let $X$ and $Y$ be the vertex sets of the component of $G'$ containing $v$ and $w$, respectively. Possibly, $X=Y$.  Graph $G'-X-Y$ has $2n-x-y$ vertices and at most \[2n-3-(n-1)-(x-1)-(y-1)=n-x-y\] edges. So, it has at least $(2n-x-y)-(n-x-y)=n$ components, and thus at least $x+y$ of them have no edges, i.e. are singletons. Either at least $x$ of them are in $V_2$ or at least $y$ of them are in $V_1$.  Suppose the former holds (the proof of the other case is symmetric). Then we place $v$ on $w$, the white neighbors of $v$ on singletons in $V_2$, and consider the remaining subgraph $G''$ with parts $G''_1$ and $G''_2$. Any packing of $G''_1$ with $G''_2$ does not create conflicts in our placement. Furthermore, $G''$ has $2n-2x$ vertices and at most $2n-3-(n-1)-(x-1)-b=n-1-x-b$ edges.  Again, the number of edges in $G''$ is less than the sizes of parts.
\qed

\begin{lem}\label{white star}
$\Delta_1,\Delta_2\leq n-3$.
\end{lem}

\noindent\textbf{Proof:}  Suppose $v,v'\in V_1$ and $N_1(v)=V_1-v-v'$.

\textbf{Case 1:} \emph{There is $w\in V_2-N(v)$ with no neighbors in $V_2$.} Send $v$ to $w$. Any packing of the resulting triple $G'=(G_1-v,G_2-w,G'_3)$ extends to a packing of $G$. Since $G'$ has at most $2n-3-(n-2)=n-1$ edges, by Corollary~\ref{n edges}, it packs unless it has a vertex of yellow degree $n-1$. But this is not the case by Lemma~\ref{yellow star}.

\textbf{Case 2:} \emph{Every $w\in V_2-N(v)$ has a white neighbor.} Let $W'$ be the set of vertices in $V_2$ reachable in $G$ from $V_1$, and let $W=V_2-W'$. Since $G-W$ has at least $(n-2)+|W'|$ edges, $|W'|\leq n-1$. So $W\neq \emptyset$ and if the white degree of $v'$ is $a$, then
\begin{equation}
\label{eq:la1}
|E(G[W])|\leq (2n-3)-(n-2)-a-|W'|=|W|-1-a.
\end{equation}

Let $W_1$ be the vertex set of a smallest tree component in $G[W]$, $y$ be a vertex of degree $1$ in $G[W_1]$ and $y'$ be the white neighbor of $y$.  Suppose the white degree of $y'$ is $b$. We send $v$ to $y$, $v'$ to $y'$ and add $a(b-1)$ yellow edges connecting the white neighbors of $v'$ with the (necessarily white) neighbors of $y'$ distinct from $y$. If the resulting triple $G'=(G_1-v-v',G_2-w-w',G_3')$ packs, then because of the added edges, this extends to a packing of $G$. Suppose it does not.  Triple $G'$ has $2(n-2)$ vertices and at most 
\begin{equation}
\label{eq:la2}
2n-3-(n-2)-a-b+a(b-1)=n-1-2a+b(a-1)
\end{equation} edges.
 
 If $a \leq 1$, then $b(a-1)\leq0$.  Also, \eqref{eq:la2} is at most $n-2$, and by Corollary~\ref{n edges}, either some $V_i$ has a vertex $z$ with yellow degree $n-2$, or the new graphs are each $K_{2}$.  However, if each of the new graphs is $K_{2}$ and $v$ originally had white degree $2$, then the case condition implies that there are at least $6 > 2n -3$ edges in the original graph.  So we assume that $z \in V_{i}$ has yellow degree $n-2$ in $G'$.  In this case, we need $a=0$ and all edges of $G$ apart from $yy'$ are incident either with $v$ or with $z$ and $z\notin \{v,v'\}$.  Thus, vertices in $V_2-y-y'$ have no white neighbors, a contradiction to the case.  Further, if $b = 1$, then there are at most $n-3$ edges in the resulting graph so $G'$ packs by Corollary~\ref{n edges}.  So let $a\geq 2$ and $b\geq 2$. In particular, $2\leq |W|\leq n$.  By~\eqref{eq:la1}, $G[W]$ has at least $a+1$ tree components, $3 \leq b + 1 \leq |W_1|\leq |W|/(a+1)\leq n/(a+1)$ and thus $2 \leq b\leq -1+n/(a+1)$. Since $a \geq 2$, then

{\allowdisplaybreaks
\begin{align*}
 |E(G')| &\leq n-1-2a+\left(\frac{n}{a+1}-1\right)(a-1) \\
 &=n-3a+n\frac{a}{a+1} - \frac{n}{a+1} \\
 &\leq n + n \frac{a}{a+1} -3a - 3 \\
 &\leq n + n\frac{a}{a+1}  - 9<  2(n-2)-3 
\end{align*}
}

Since  $G'$ does not pack, by induction (the last strict inequality ensures that examples from Figure~\ref{fig:B-E} do not appear) some vertex $z$ in $G'$ has $d_3(z)=n-2$ or $d_i(z)=n-3$ for some $i=1$ or $2.$ But $d_i(z)\neq n-3$, since we deleted at least $n-2+a+b\geq n+1$ edges from $2n-3$ in $G$ and have not added white edges.  Similarly, since we have not added any yellow edges incident to $V_{2} - W_{1}$, we cannot have $d_3(z)=n-2$ if $z\in V_{2} - W_{1}$.

The case when $z\in V_1$ and has yellow degree $n-2$ is also forbidden since for this to happen, every vertex in $V_{2}$ must be incident to a yellow edge.  However, $G[W]$ has at least $a+1\geq 3$ components, of which only $W_{1}$ is incident to yellow edges.  Finally, we cannot have  $z\in W_{1}$, or else we must have created a star by adding yellow edges from $z \in W_{1}$ to all vertices in $V_{1} - v - v'$.  But this implies a packing by Theorem~\ref{B-E} since if $a = n -2$,  then $e_{1} = 2n-4$ and $e_{2} = 1$ and $G$ contained no yellow edges. \qed

\begin{lem}\label{white neighbor}
Every vertex of $G$ has a white neighbor.
\end{lem}

\noindent \textbf{Proof:} Suppose $v\in V$ has no white neighbor.
	
\textbf{Case 1:} \emph{$v$ is isolated in $G$.} Without loss of generality, assume $v\in V_1.$ If any $w\in V_2$ has degree at least $2$ in $G$ then placing $v$ on $w$ decreases $e_1+e_2+e_3$ by at least $2$. If the new triple packs, then this extends to $G$, otherwise by Lemmas~\ref{yellow star} and~\ref{white star}, it is one of the examples from Figure~\ref{fig:B-E}.  Suppose that $(G_{1} - v, G_{2} - w)$ is an example from Figure~\ref{fig:B-E}, then there must be exactly $2(n-1)-3$ edges and so $d(w) = 2$.  If both edges adjacent to $w$ were white, then $G$  originally had no yellow edges and the theorem holds by Theorem~\ref{B-E}.  If both edges were yellow, since every graph in the Figure~\ref{fig:B-E} contains 3 vertices of positive degree, there is some $v' \in V_{1} - N(w)$ with degree at least $1$.  Then $(G_{1}-v', G_{2} - w)$ contains strictly fewer than $2(n-1)-3$ edges and therefore, by Theorem~\ref{B-E}, pack.  This packing can be extended to a packing on $G$.  Finally assume that $w$ has exactly one neighbor $w' \in V_{2}$ and one neighbor in $V_{1}$.  Again, we can choose some $v' \in V_{1} - N(w)$ with positive degree.  Create a new graph triple by removing $v'$ and $w$ and adding yellow edges from $w'$ to $N(v')$.  This new triple has exactly $2(n-1)-3$ edges.  By Lemmas~\ref{yellow star} and~\ref{white star} (and since it has at least one yellow edge), the new triple packs and can be extended to a packing of $G$.

So suppose the total degree of each $w\in V_2$ is at most $1$. If at least one $w\in V_2$ is isolated, then symmetrically, each $v'\in V_1$ has degree at most $1$ and $G$ packs by Corollary~\ref{n edges}. Thus $d(w)=1$ for each $w\in V_2$. If there are no yellow edges, then we are done by Theorem~\ref{B-E}. So let $wv'\in E_3$. If all vertices in $V_1-v'$ are isolated, then the total degree sum of $G$ is at most $n+\Delta_3(G)\leq 2n-3$ and so $e_1+e_2+e_3<n$, a contradiction to Corollary~\ref{n edges}. Otherwise, let $u \in V_1 - v'$ be a vertex of maximum degree and send $u$ to $w$.  If $(G_1 - u, G_2 - w)$ packs, then this packing extends to a packing of $G$.  If it does not pack, then by induction, $(G_1 - u, G_2 - w)$ is an example from Figure~\ref{fig:B-E} and $d(u) = 1$.  However, each example in Figure~\ref{fig:B-E} contains a graph with multiple vertices of degree at least 2, contradicting the maximality of $u$.

\textbf{Case 2:} \emph{$N_{3-i}(v)\neq\emptyset$.} Let $v \in V_1$ such that $N_1(v) = \emptyset$ and suppose $w' \in N_2(v)$.  Among the vertices in $V_2 - N_2(v)$ with maximal degree, let $w$ be a vertex that minimizes $d_3(w)$.  We send $v$ to $w$ and consider the triple formed by removing these two vertices.  If $d(v) + d(w) > 2$, then the remaining graph triple packs by induction and the packing extends to a packing of $G$.  Therefore, by Case 1, we may assume that $d(v) = d(w) = 1$.   By induction, $(G_1 - v, G_2 - w)$ must be an example from Figure~\ref{fig:B-E}, as otherwise the graphs pack.

However, by the maximality of $d(w)$, all vertices in $V_2 - w'$ must have degree at most $1$ in $G$ and, hence, in $G_2 -w$.  By inspection, $G_2 - w$ is either G(1) or G(3) in Figure~\ref{fig:B-E}, as all other graphs have multiple vertices with degree at least $2$.  Since H(1) and H(2) each have an isolated vertex that, by Case 1, was not isolated in $G_1$, we must have removed an incident yellow edge when deleting $w$.  Both G(1) and G(3) have at least $4$ vertices adjacent to at exactly one white edge.  In the process of removing $v$ and $w$ from $G$, we have removed at most $2$ edges incident to $V_{2}$.  Thus, in $G_2$, there must also have been a vertex with degree $1$ adjacent to a white edge, contradicting our choice of $w$.
\qed

\noindent \textbf{Proof of Theorem~\ref{List B-E}:}
Let $G$ be our minimum counterexample. If $G$ has no yellow edges, then the original Theorem~\ref{B-E}~applies. So suppose $G$ has a yellow edge  $xy$ with $x\in V_1$ and $y\in V_2$. Since $|E(G)|\leq 2n-3<2n$, there are vertices of degree at most $1$.  We may assume that $v\in V_1$ and $d(v)\leq 1$.  By   Lemma~\ref{white neighbor}, $v$ has a white  neighbor, $v'$ (possibly $v'=x$). We send $v$ to  $y$ and add yellow edges from $v'$ to
each white neighbor of $y$. Then we obtain the triple $G'$ with exactly two edges less.  Since we have at least one yellow edge (connecting $v'$ with a white neighbor of $y$), we do not get a graph from Figure~\ref{fig:B-E}.  So by Lemmas~\ref{yellow star} and~\ref{white star}, the theorem is proved. \qed

%
%
%
%

\bigskip

{\bf Acknowledgement:} We would like to thank Gexin Yu for helpful comments.

\bibliographystyle{abbrv}
\bibliography{ListPacking}

\end{document}